\magnification=\magstep1
\overfullrule=0pt 

\centerline{\bf Measuring Large Scale Space Perception in Literary Texts}
\bigskip
\bigskip
Paolo Rossi
\smallskip
{\it University of Pisa, Pisa, Italy}
\bigskip
{\it Address for correspondence}: [P.Rossi, Dipartimento di Fisica ``E.Fermi'', L.Pontecorvo 3, 56127 Pisa, Italy]

E-mail: Paolo.Rossi@df.unipi.it
\bigskip
{\bf Summary.} The center and radius of perception associated with a written text are defined, and algorithms for their computation are presented. Indicators for anisotropy in large scale space perception are introduced.
The relevance of these notions for the analysis of literary and historical records is briefly discussed and illustrated with an example taken from medieval historiography.
\bigskip 
\bigskip
{\bf 1. Introduction}
\smallskip
All attempts to represent human perceptions by quantitative and measurable parameters are certainly quite difficult as well as conceptually debatable.

However, when working on such cultural artifacts as written records, a partial simplification occurs because of the substantially limited and selected amount of information encoded in the data.
Therefore trying to evaluate parameters from written records is like working in a very controlled laboratory environment.

Quantitative study of written language has a long history, but we are not aware of any attempt to define quantitative indicators of subjective perception, that might be measured on any given record.

We think that an author's (conscious or inconscious) perception of space might be evaluated by a statistical study of topographic and/or geographic occurrences in the text. In particular we hope to be able to give convincing quantitative definitions of the ``center of perception'' and ``radius of perception'' for a given text, measuring respectively the author's vantage point and the typical distance scale involved in any specific document.

Both quantities can certainly be qualitatively estimated, in most instances, by a critical reading of the document, but we believe that a formal algorithm may help to remove all possible contamination caused by the influence of the author's style or even by the reader's own radius of perception, and may therefore provide valuable historical and/or sociological insights.

We therefore propose that such parameters (and a few auxiliary ones we shall define later) might be of some use in the study of a broad range of written records, from fiction to historiography. In particular, when applied to historiography, they might be helpful in classifying and ranking primary sources of historical information, and in evaluating the social and time evolution of the perception of geographical scales.
\bigskip
\eject
{\bf 2. Data organization}
\smallskip
Two preliminary actions must be performed on the records for the application of our algorithm: 

1) The record must be indexed, and all occurrences of topographic and/or geographic names must be classified and counted. In the simplest and most effective version of our indicators, only names that can be associated with point-like entities (towns, villages, buildings, individual mountains, ...) must be retained, while extended geographical objects (states, regions, rivers, mountain chains, ...) do not easily lend themselves to simple quantitative manipulation.

2) Coordinates of individual geographic entities must be identified. 

We note here that in practice most nonfiction books are already (more or less carefully) indexed, and that very large online databases of geographic coordinates are available. This preliminary data organization should therefore usually be available with comparatively small effort.
\bigskip
{\bf 3. Definitions and basic notions}
\smallskip
We first introduce our notation. Let $o_n$ be the number of occurrences of the n-th geographic item appearing in a given record, $\theta_n$ its latitude and $\varphi_n$ its longitude.

It is also convenient to define the ``weight'' $w_n$ of the n-th item by the standard relationship
$$w_n = {o_n \over \sum_n o_n}.$$

We first identify the ``center of perception'' of a given record. In most cases, this is more or less explicitly declared by the author, or it is easily detectable by qualitative inspection of the text.
It is, however, convenient to give a formal definition, which may  be contrasted with qualitative information, in order to check consistency.
Defining  the weighted coordinates as:
$$X = \sum_n w_n \cos \theta_n \cos \varphi_n,\qquad Y = \sum_n w_n \cos \theta_n \sin \varphi_n,\qquad Z = \sum_n w_n \sin \theta_n,$$
the geographic coordinates of the center of perception can now be defined as
$$\Theta = \arctan {Z \over \sqrt{X^2+Y^2}},\qquad \qquad \qquad \Phi = \arctan {Y \over X}.$$

Notice that the above definition is a slight generalization of the notion of ``center of mass'' of a physical system, adapted in order to account for the two-dimensional and spherical nature of the Earth's surface.
The distance $d_n$ of any given point from the center of perception is then defined by:
$$d_n = R_E \arccos [\cos \theta_n \cos \Theta \cos (\varphi_n - \Phi) + \sin \theta_n \sin \Theta ],$$
where $R_E$ is the Earth's radius. When only small portions of the Earth's surface are involved, it is quite easy to derive approximate expressions for $\Theta$, $\Phi$ and $d_n$, but in practice the use of approximate expressions is unnecessary, since the exact expressions can easily be evaluated by the use of rather standard computer programs (worksheets).
\smallskip
We are now ready to define the ``radius of perception'', $R$, as the weighted average of the distances:
$$R = \sum_n w_n d_n.$$

This definition can obviously be applied not only to the ``mathematical'' center of perception but also to the ``empirical'' or intuitive one. If there is no substantial discrepancy between the two locations, the resulting values of the radius of perception are expected to agree closely.

In order to get a more direct interpretation of our definition, let's notice that, for most nonsingular distributions, $R$ is the radius within which about half of the geographic quotations present in the text can be found. More precisely, for an isotropic Gaussian distribution, the exact fraction of quotations included within the radius of perception is
$$1 - e^{-{\pi \over 4}} = 0.544...$$

Inspired by the standard analysis of mass distributions, we can also find a quantitative description of the anisotropy (orientation dependence) in space perception by defining the ``ellipse of perception'' and its principal axes (``axes of perception''), corresponding to a rather straightforward generalization of the inertia tensor and its properties.
We introduce, as an intermediate step in the computation, a new set of spherical coordinates $\theta'_n$, $\varphi'_n$, such that the center of perception turns out to be the North Pole of the new coordinate system. The new coordinates are obtained from:
$$\theta'_n = {d_n \over R_E},\qquad \varphi'_n = \arctan \Bigl[{\cos \theta_n \sin  (\varphi_n-\Phi) \over \cos \theta_n \sin \Theta \cos (\varphi_n-\Phi)- \sin \theta_n \cos \Theta} \Bigr].$$

Once the transformation has been performed, one must reinterpret $d_n$ and $\varphi'_n$ as {\it planar} polar coordinates, and construct a symmetric tensor whose components are
$$I_{11} \equiv I_+ + I_-, \qquad I_{22} \equiv  I_+ - I_-, \qquad I_{12} = I_{21},$$
where we have defined
$$I_+ = {1 \over 2}\sum_n w_n d_n^2, \qquad I_- = {1 \over 2}\sum_n w_n d_n^2 \cos 2 \varphi'_n, \qquad  I_{12} = {1 \over 2}\sum_n w_n d_n^2 \sin 2 \varphi'_n. $$
Notice that $I_+$ and $I_-^2 + I_{12}^2$ are invariant under rotations. 

The axes of perception are two orthogonal directions, respectively forming with the NS and EW geographic axes an angle $\Phi'$ identified by the relationship
$$\Phi' = {1 \over2} \arctan {I_{12} \over I_-}.$$

The eigenvalues of the tensor are the square lengths of the two semiaxes of the ellipse. Their (rotation invariant) values are
$$a^2 = I_+ + \sqrt{I_-^2+I_{12}^2}, \qquad \qquad \qquad b^2 = I_+ - \sqrt{I_-^2+I_{12}^2}.$$
The lengths $a$ and $b$ are obviously related to the radius of perception, but the relationship depends on the details of the distribution; in the simple case of Gaussian distributions with slight anisotropy the following relation holds:
$$(a^2+ b^2) = {4 \over \pi} R^2.$$

The ellipse of perception offers further significant information on the geographic perception of the author, since it can reveal any privileged direction in the focus of attention and/or narration.
This will be especially true when the lengths of the two axes turn out to be significantly different.
\bigskip
{\bf 4. Treatment of extended regions}
\smallskip
Our definitions rely heavily on the notion of geographic coordinates, and cannot be trivially applied to extended geographical objects. 
However, since in many cases the number of quotations of extended regions is numerically comparable with the number of references to pointlike geographical objects, it is certainly worth trying to exploit such a relevant source of information.

The treatment of extended objects should be performed independently of the analysis of pointlike ones. The two results can eventually be compared to check consistency, or to ferret out the sources of inconsistency.

The first and simplest possibility involves assigning to each extended object the coordinates of its geometrical center. This procedure is formally consistent, but in practice it only makes sense when all objects involved have comparable extensions. Moreover, one must keep in mind that the resulting values of the estimated parameters are necessarily affected by a statistical error whose magnitude can be roughly evaluated by taking the ratio between the average radius of the regions in the sample and the square root of their number.
With these warnings, it is certainly useful to check if the figures thus obtained, which are certainly independent of those extracted from the analysis of pointlike objects, are in reasonable agreement with them.
In this context we must also mention the possibility of more subjective, but sometimes more significant, approaches. It is indeed possible to associate to each extended region the coordinates of places which may have been especially significant for the author, independent of their ``centrality'' (just think of the role of political centers for states, mountain passes for mountain chains, main or nearest bridges for rivers, and so on).
\bigskip
{\bf 5. A question of metric and measures}
\smallskip
We are obviously aware of the fact that distance perception in individuals is strongly conditioned by several psychological and social factors, implying not only a substantial anisotropy, but also an ``effective (and subjective) metric'' which is certainly not Euclidean; the ratio between the perceived extension of a region and the real one usually tends to decrease with distance, and this is only the best known example.

These considerations should not, however, be seen as a real obstacle against adoption of our formalized definitions. Indeed the use of geometric coordinates and standard measure units should not obscure the fact that the notions presented here are essentially of a topological nature. We already stressed that the radius of perception identifies the region including (roughly) half of the geographic quotations. Expressing the radius in units that are familiar to us allows easy comparison between texts, but {\it within} a given text it would be more proper to treat the radius as an independent measure unit, defining the scale for an intrinsic evaluation of other distances within the same record.

Another approach to the identification of appropriate distance scales could be based on the study of contemporary sources containing indications of perceived distance, especially travel reports including a specification of the times involved in diplacements.

Concerning anisotropies, we think that the determination of the ellipse of perception can give a first substantial contribution to their identification. In principle, more accurate information might be extracted from the data by applying a principal component analysis. This however would require sufficiently large samples in order to be statistically reliable.
\bigskip
{\bf 6. Text correlation}
\smallskip
When comparing two authors or even two different texts by the same author it is possible to use the previous definitions to perform quantitative comparisons. But it is also possible, in the cases of interest, to express the results of comparison in terms of another quantitative indicator, which can be defined as ``text correlation''.

Let $\{p_n\}$ and $\{q_n\}$ be the sets of occurrences (or weights) of two different records. Then, following the usual statistical meaning, text correlation $C_{pq}$ is defined by
$$C_{pq} = {\sum_n p_n q_n \over \sqrt{\sum_n p_n^2} \sqrt{\sum_n q_n^2}}.$$

By construction $C_{pq}$ is a number ranging from 0 to 1, where 0 represents total decorrelation (no term quoted in a text is quoted in the other), while 1 represents full correlation (every term is quoted with exactly the same relative frequency in both texts). Note that full correlation by no means implies identity of structure and content: completely different sentences can be made out of the same words used with the same frequency.
In practice it is convenient to introduce also the ``text decorrelation'', whose self-explanatory definition is $D_{pq} \equiv 1 - C_{pq}$.
The notion of text correlation (or decorrelation) is completely general, and it might therefore be very useful well beyond the specific domain of the present discussion, since the analysis may be extended to all the terms appearing in a text, or at least to all the ``strongly significant'' terms (nouns, adjectives, verbs and adverbs).

In the context of our space perception analysis, text correlation may be employed in the analysis of texts having some (known or possible) genealogical link.
A strong correlation may be used as a confirmation of the link, but also, subtracting from one text the correlated component with standard techniques of vector calculus, one may focus on the existing differences and try to trace their origin. 
\bigskip
\eject
{\bf 7. An illustrative example}
\smallskip
The history of Northern France in the tenth century has been reported to us essentially by only two contemporary writers, both monks and both living most of their lives in the city of Rheims: the older is Flodoard, whose {\it Annals} cover the years 919-966, and the younger is Richer, whose {\it Histories}, divided into four books, go from 888 to 998.

Richer declaredly borrows information from Flodoard, but since the timespans covered by the two are significantly different, it may be of some historical and psychological interest to compare the values of the above defined indicators in the two texts.

First of all, it is easy to verify that the ``mathematical'' center of perception of both authors, evaluated according to our definition, is located in the immediate vicinity of Rheims. This may appear as a trivial statement for any actual reader of the texts, since the ``intuitive'' result is absolutely self-evident, but this is certainly a first test of accuracy, which offers also an indication for an estimate of the systematic error.

We evaluated the radius of perception, obtaining the following results:

Flodoard's radius: 178 Km (based on 545 quotations)

Richer's radius: 201 Km (based on 658 quotations)
\smallskip
We also evaluated the radius of perception for (roughly 500 quotations of) extended regions. Our evaluation was very rough, since for distances of geographic areas from Rheims we kept only the first significant digit, discounting statistical errors in the data not smaller than 50 Km. Since about thirty regions were involved, we expected our final error to be as big as 10 Km. With these specifications, it was encouraging to find that Flodoards's regional radius of perception is about 210 Km and Richers's is about 200 Km.

We then performed the construction of the ellipse of perception, and found  that the main axis of perception connects directly Rheims with Rome. More precisely, Flodoard's value of $\Phi'$ is $45^\circ$, while Richer's value is $44^\circ$, and Rome's orientation with respect to Rheims is $42^\circ$. 

Again this may appear obvious, since we are analyzing medieval authors who are monks and are very interested in all events related to religious life.
However this observation, together with the numerical evidence for a singularity in the distributions related to quotations of Rome, suggests that we can remove all reference to Rome from both lists (more than 20 quotations in each text) and evaluate a ``biased'' radius of perception referred to the residual data. Note that this is an admittedly naive (but not too inaccurate) way of performing the first step of a principal component analysis.

Interestingly enough, the resulting values are dramatically lower:

Flodoard's biased radius: 134 Km

Richer's biased radius: 154 Km
\smallskip
The new ellipses of perception are completely reoriented, and they now show a preferred direction which is no longer the same for the two authors, in both cases connecting Rheims, in a NE-SW alignment, to some of the most important centers of political power, but with an apparent shift in the authors' focus (on which we shall comment later). The new values of $\Phi'$ are $-85^\circ$ for Flodoard (with axes lengths $a$ =  156 Km and $b$ = 135 Km) and $-61^\circ$ for Richer (with  $a$ = 167 Km and $b$ = 152 Km). 

In passing we notice that, in the light of these results, the celebrated ``journey to Chartres'' reported in Richer IV, 50 should be viewed in a slightly different light, keeping in mind that the distance between Chartres and Rheims is about 200 Km, a value which locates the endpoint of Richer's travel at the border of his perceptive space.

Concerning the question of ``intrinsic'' distance scales, we may find some help in the study of a contemporary document: the detailed report of a journey from Rome to Canterbury through Rheims, made by archbishop Sigeric sometime between990 and 994 along the {\it via Francigena} and involving eighty stations. Eight of these stations are quoted also by Flodoard and Richer. The average distance between two stations (North of the Alps) is about 27 Km, and therefore it is meaningful to say that the radius of perception of the two authors under scrutiny corresponds to the distance that could be reached in the tenth century in a 5 or 6 days' walk.
\smallskip
The (geographic) text correlation between Flodoard and Richer is an impressive 0.96, indicating that there is an almost total superposition of the geographic name distributions.
This makes any difference between the two data sets especially significant.

In particular, subtracting from the set of Richer's geographic occurrences the projection of Flodoard's, we are left with two small clusters concentrated in two quite restricted areas: Normandy (7 occurrences) and the region around Li\`ege (10 occurrences).
Both areas are at the borders of Richer's perceptive space, and therefore these occurrences must have some structural explanation. While the political role of the Normans in tenth century France may account for the first set, the focus on the region around Li\`ege is better explained, with the help of other circumstantial evidence which we shall not discuss here, by the conjecture that this may be Richer's own ancestral area.

On the other hand, by subtracting from Flodoard's occurrences the projection of Richer's, we found about fifty extra references to places in the Rheims region, which Richer doesn't mention. The link of Flodoard to Rheims area, whether strictly personal or ancestral, is certainly much more marked than that of Richer.

 Without belaboring the point, we only want to stress that these facts seem to have escaped all previous (qualitative) analysis and comparison of these records.
\smallskip
We would also like to make a short comment on the results of a broader application of the notion of text correlation, which was evaluated the lists of all words, and of significant words, extracted from the four books of Richer's.

The  decorrelation coefficients for all words are
$$D_{12} = 0.066 \qquad \qquad D_{34} =  0.068$$
$$D_{13} = 0.125 \qquad D_{14} = 0.153 \qquad D_{23} = 0.103 \qquad D_{24} = 0.131$$
while the coefficients for significant words are
$$D_{12} = 0.136 \qquad \qquad D_{34} = 0.153$$
$$ D_{13} = 0.212 \qquad D_{14} = 0.184\qquad D_{23} = 0.184 \qquad D_{24} = 0.174$$

Correlation between Book 1 and 2 and correlation between Book 3 and 4 are quite strong, but the first and second couple of books are significantly decorrelated. This quantitative hint might suggest some deeper qualitative analysis of the texts and their context.
\bigskip
{\bf 8. Conclusions}
\smallskip
Treating geographic occurrences in literary texts with purely formal and statistical techniques might appear as an attempt to give oversimplifying quantitative answers to very subtle and essentially qualitative questions.

Nevertheless, we hope to have been able to present sufficiently simple algorithms and sufficiently convincing evidence for the opportunity of starting a systematic exploration based on the measurement of the indicators we have proposed.

Furthermore, a detailed quantitative study of the distribution of geographic occurrences, especially when performed on sufficiently large samples, might help to construct, for a given author or epoch, a hierarchy of ``perceived distances'' that might turn out to be significantly different from that resulting just from measurement of physical distances.
The role of Rome and the focus on different regions by different authors, which we discussed in our case study, are just specific instances of this notion of ``perceived distance'', and show some of the ways in which our methods may be helpful in identifying this phenomenon.

Finally, we would like to mention that the extension of some of the ideas discussed here to time perception and time distance is conceptually staightforward and mathematically simple. Practical problems may arise for lack of ``chronological indexing'' of the texts and because of the often overwhelming presence of ``extended'' instead of ``pointlike'' time references. These difficulties should not, however, prevent attempts in this direction.
\bigskip
{\bf Acknowledgments}
\smallskip
I am deeply indebted with Steve Shore for his critical reading of the manuscript, but especially for his many subtle comments and stimulating suggestions, which triggered some important extensions and revisions. I am also indebted with Gianfranco Prini, whose skeptical attitude forced me to find answers to previously overlooked but crucial questions.
\bigskip
{\bf References}
\smallskip
Flodoard (1906) {\it Annales}, ed. by P.Lauer as vol. XXXIX of {\it Collection de textes pour servir a l'enseignement de l'histoire}

Richer (2000) {\it Historiae}, ed. by H.Hoffmann in {\it M.G.H. Scriptores} 38

Sigeric (1874) {\it Memoria}, ed. by W.Stubbs in {\it Rerum Britannicarum Medii Aevi Scriptores}, vol. 63, pp. 391-399
\bigskip
\bigskip
\bigskip
{\bf Figure caption}

The map covers a wide area in Western Europe, between the rivers Loire and Rhine, and reports geographical names explicitly quoted by Richer. The circles and straight lines represent Richer's (unbiased and biased) perception radii, and Richer's (biased) principal axes of perception.
The continuous and broken lines represent repectively the approximate track of the {\it Via Francigena} (with Sigeric's stations) and the path of Richer's journey to Chartres.
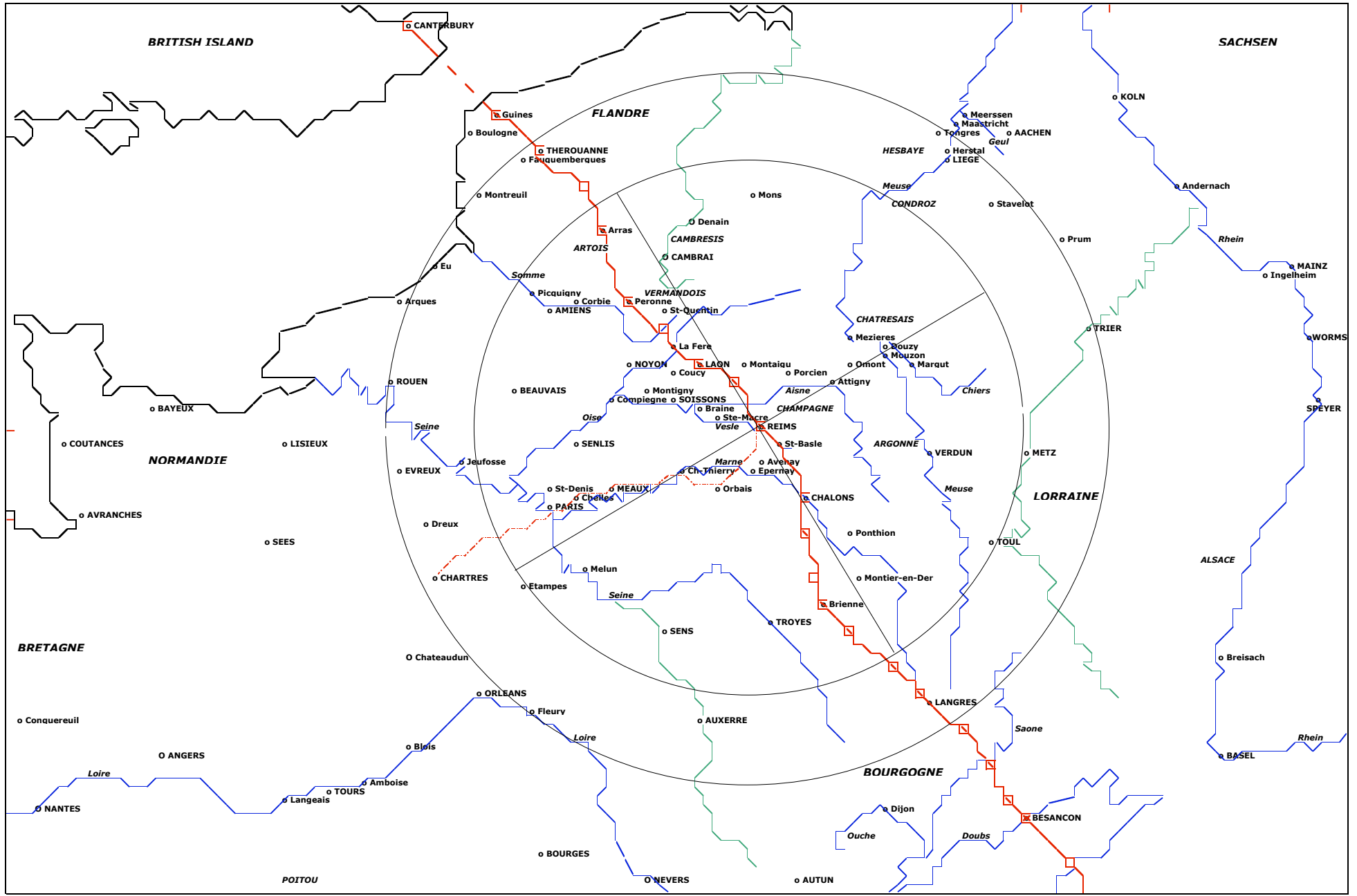
\end